# *BourbOulipo*
# *Relazioni tra Oulipo e Bourbaki*

di

Elena Toscano e M. Alessandra Vaccaro

Dipartimento di Matematica e Informatica, Università degli Studi di Palermo

> «La struttura è libertà»
> *I. Calvino*

**Abstract**

The theme of the influence of Bourbaki's ideas and methodologies on disciplines that transcend the mathematical field has been frequently addressed and, in particular, the reflection on the relationship between Bourbaki and Oulipo has animated a passionate debate within literary criticism.

The purpose of this contribution is to investigate the relations between Le Lionnais, Queneau and Roubaud - and more generally Oulipo - and Bourbaki. Through a comparative examination between the two groups, the existence of an undeniable charisma of the "polycephalic mathematician" on the artistic-literary movement was highlighted.

**Abstract**

Il tema dell'influenza delle idee e delle metodologie bourbakiste su discipline che trascendono l'ambito matematico è stato frequentemente affrontato e, in particolare, la riflessione sul rapporto tra Bourbaki e Oulipo ha animato un appassionato dibattito in seno alla critica letteraria.

Lo scopo di questo contributo è quello di investigare le relazioni tra Le Lionnais, Queneau e Roubaud - e più in generale l'Oulipo - e Bourbaki. Attraverso un esame comparativo tra i due gruppi si è evidenziata l'esistenza di un innegabile carisma del "matematico policefalo" sul movimento artistico-letterario.

## Introduzione

In [42] abbiamo voluto delineare la figura di François Le Lionnais (Paris, 3 ottobre 1901 – Boulogne-Billancourt 13 marzo 1984) mediante l'analisi di alcuni dei suoi contributi più significativi che riflettono la sua concezione del profondo legame esistente tra matematica e letteratura. Le Lionnais, ingegnere chimico di formazione, durante gli anni giovanili si avvicina al Dadaismo e diviene amico di Marcel Duchamp con cui condivide la passione per gli scacchi. Membro della Resistenza, durante la seconda guerra mondiale è deportato nel campo di concentramento di Dora. Gli anni successivi alla fine della guerra lo vedono intento nella concretizzazione de *Les grands courants de la pensée mathématique* [17] un'opera sullo stato dell'arte, gli sviluppi delle ricerche in matematica e le sue influenze in altri ambiti del sapere, il cui progetto - stimolato dall'amico editore Jean Ballard - era stato avviato già durante il periodo dell'occupazione nazista[1].

Le Lionnais riesce nell'intento non banale di coinvolgere e coordinare non solo alcuni tra i più grandi matematici dell'epoca, tra cui esponenti di Bourbaki quali André Weil, Jean Dieudonné e Roger Godement, ma anche artisti e intellettuali di spicco come Le Corbusier e Louis De Broglie.

La preparazione de *Les grands courants de la pensée mathématique* è l'occasione per entrare in contatto nel 1943 con Raymond Queneau (Le Havre, 21 febbraio 1903 – Parigi, 25 ottobre 1976) col quale stringerà un'amicizia e un sodalizio intellettuale che dureranno per tutta la vita e che, nel 1960, vedranno la loro massima espressione nella fondazione dell'Oulipo (*Ouvroir de littérature potentielle*). L'idea di una possibile sperimentazione matematica nel processo creativo letterario è concepita da Le Lionnais nel corso dei suoi studi universitari e, dopo averlo continuamente affascinato per molti anni, poté realizzarsi grazie alle conversazioni con Queneau. Lo scopo principale dell'*Ouvroir* è investigare l'applicazione di alcune strutture della matematica ai diversi livelli della produzione letteraria, grazie all'imposizione di vincoli e regole (*contrainte*) che, lungi dall'imbrigliare l'estro creativo degli scrittori, permettono di dare origine a nuove forme espressive. Le Lionnais è in primo luogo il 'pensatore' dell'Oulipo dando dei contributi specifici soprattutto nella stesura di tre *Manifesti*[2] e in una serie di lavori pubblicati nell'opera collettiva *La littérature potentielle* [28]. Insieme a Queneau, Le Lionnais gioca anche il ruolo importante di reclutatore dei cosiddetti *oulipiens* al cui primo nucleo di dieci soggetti si uniscono, alcuni anni dopo, anche gli scrittori Georges Perec, Italo Calvino e Jacques Roubaud.

---

[1] Per una biografia dettagliata di Le Lionnais si rimanda, per esempio, a [39], [40] e [41].

Negli anni in cui nasce l'Oulipo il gruppo di Bourbaki, che è all'apice del successo, aveva creato una vera e propria rivoluzione nel modo di concepire la matematica. Infatti in quegli anni Bourbaki è una realtà già consolidata che, con il pretesto di rinnovare e riorganizzare l'insegnamento superiore della Matematica, aveva finito per dare un nuovo impulso alla staticità dell'ambiente matematico francese. Nel 1950 Le Lionnais puntualizza «quando il movimento bourbakista è comparso, la matematica rischiava di soccombere sotto il disordine della propria ricchezza. Secondo gli autori - da un paese all'altro e in uno stesso paese - le stesse idee venivano espresse mediante parole e simboli differenti e le stesse parole o gli stessi simboli esprimevano delle idee differenti. La storia delle scienze riconoscerà senza dubbio che uno dei principali meriti di Bourbaki è d'avere coraggiosamente iniziato a mettere fine a quest'anarchia fissando la terminologia matematica»[3].

Bourbaki contribuisce a *Les grands courants* con un articolo che si può considerare il chiaro manifesto ideologico del gruppo, come si evince dall'affermazione «l'evoluzione interna della matematica, malgrado le apparenze, ha rafforzato più che mai l'unità delle sue diverse parti e ha creato una sorta di nocciolo centrale più coerente di quanto sia mai stato. L'essenziale di questa evoluzione è consistito in una sistematizzazione delle relazioni esistenti tra le diverse teorie matematiche e si riassume in una tendenza che è generalmente conosciuta sotto il nome di "metodo assiomatico". […] È solamente con questo senso della parola "forma" che si può dire che il metodo assiomatico è un "formalismo"; l'unità che essa conferisce alla matematica non è l'armatura della logica formale, unità di scheletro senza vita; è la linfa che nutre un organismo in pieno sviluppo, il duttile e fecondo strumento di ricerche al quale hanno coscientemente lavorato, dopo Gauss, tutti i grandi pensatori della matematica, tutti coloro che, seguendo l'espressione di Lejeune-Dirichlet, hanno sempre cercato di "sostituire le idee al calcolo"».[4]

Il tema dell'influenza delle idee e delle metodologie bourbakiste su discipline che trascendono l'ambito matematico è stato frequentemente affrontato e, in particolare, la riflessione sul rapporto tra Bourbaki e Oulipo ha animato un appassionato dibattito in seno alla critica letteraria. Nonostante tale relazione sia addirittura suggerita da alcuni *oulipiens* quali Le Lionnais e Roubaud in testa oltre che Marcel Bénabou [13], l'argomento probabilmente è stato spesso trattato in modo incompleto e soprattutto semplicistico suggerendo *tout court* che vi sia una "relazione di subordine" di Oulipo nei confronti di Bourbaki evidentemente basata anche su motivazioni meramente cronologiche. Ciò ha indotto presumibilmente

---

[2] *Le Premier Manifeste* e *Le Second Manifeste* sono pubblicati in [28]. *Prolégomènes à toute littérature future* noto come *Le Troisième Manifeste* è stato ritrovato tra i documenti di Le Lionnais e pubblicato postumo in [26].

[3] «Lorsque le mouvement bourbakiste est apparu, les mathématiques risquaient de succomber sous le désordre de leurs richesses. Selon les auteurs - d'un pays'à l'autre et dans un même pays - les mêmes idées s'exprimaient par des mots et par des symbols différents, et les mêmes mots ou les mêmes symbols exprimaient des idées différentes. L'histoire des sciences recoinnaîtra sans doute comme l'un des principaux mérites de Bourbaki d'avoir courageusement entrepris de metre fin à cette anarchie en fixant la terminologie mathématique» [18, p.10].

[4] «l'évolution interne de la science mathématique a, malgré les apparences, resserré plus que jamais l'unité de ses diverses parties, et y a créé une sorte de noyau central plus coherent qu'il n'a jamais été. L'essentiel de cette évolution a consisté en une systématisation des relations existant entre les diverses theories mathématiques, et se résume en une tendance qui est généralement connue sous le nom de "méthode axiomatique". [...] C'est seulement avec ce sens du mot "forme" qu'on peut dire que la méthode axiomatique est un "formalisme"; l'unité qu'elle confère à la mathématique, ce n'est pas l'armature de la logique formelle, unité de squelette sans vie; c'est la sève nourricière d'un organisme en plein développement, le souple et fécond instrument de recherches auquel ont consciemment travaillé, depuis Gauss, tous les grands penseurs des mathématiques, tous ceux qui, suivant la formule de Lejeune-Dirichlet, ont toujours tendu à "substituer les idées au calcul"» [17, p. 37, p. 47].

l'*oulipienne* Michèle Audin ad esaminare e documentare il rapporto tra Oulipo e Bourbaki [4] concludendo che su di esso, il più delle volte, si pone un'enfasi eccessiva e non correttamente fondata.
Uno studio comparativo tra i due gruppi è tuttavia possibile e a esso è dedicato il presente lavoro.

## Le origini dell'Oulipo: dall'opposizione al Surrealismo all'adesione al Bourbakismo

Come si è già detto, quando Le Lionnais incontra Queneau è subito sedotto dalle sue idee[5] tanto da proporgli di collaborare alla realizzazione di *Les grands courants* e, come è noto [42], il suo contributo si concretizzerà nel lavoro *La place des Mathématiques dans la classification des Sciences*, uno tra i diciotto articoli della sezione *Influenze* [17].

Queneau - scrittore, poeta e drammaturgo - è ormai in quegli anni una figura di spicco del panorama letterario francese e internazionale. Dopo aver mosso i primi passi, da giovane studente della Sorbona, tra le fila del Surrealismo[6], affascinato dall'approccio surrealista alla creazione che si basa sull'esaltazione dell'inconscio e del subconscio, se ne allontana bruscamente (anche in seguito a dissapori personali con André Breton) tracciando una propria via ispirata alla sperimentazione matematica in letteratura.

Lo stesso Oulipo è stato fondato, nella sua struttura e nel suo funzionamento, molto chiaramente in opposizione al gruppo surrealista. Secondo una mirabile analisi [34] di Roubaud infatti la storia della letteratura francese, almeno dal Rinascimento in poi, è stata caratterizzata dall'apparire di correnti, movimenti o, più genericamente, gruppi letterari accomunati dalle seguenti caratteristiche principali:

i. la volontà di rinnovare e rifondare una letteratura ormai deteriorata dai predecessori;
ii. la divisione interna del gruppo in ranghi: da un lato uno o più leader, dall'altro i membri subordinati;
iii. il disprezzo per gli altri gruppi contemporanei;
iv. un modo di lavorare che attraverso scissioni, divergenze ed espulsioni conduce alla dissoluzione abbastanza rapida del gruppo stesso.

È dunque evidente come il gruppo surrealista, seppur involontariamente, abbia svolto il ruolo di antagonista dell'Oulipo. Le Lionnais e Queneau, nel fondarlo, non ebbero dubbi: un gruppo aperto che si accresce grazie alla cooptazione unanime di nuovi membri, nessuna espulsione o dimissione, nessuna chiusura verso le influenze culturali antecedenti e/o contemporanee.

Circa un decennio dopo la nascita del Surrealismo la Francia vedeva emergere, questa volta in ambito matematico, un altro gruppo di avanguardia noto come Bourbaki.

Nicolas Bourbaki è l'eteronimo dietro cui, a partire dal 1935, si cela un gruppo di brillanti matematici prevalentemente francesi ed ex allievi dell'*École Normale Supérieure* «che hanno deciso di realizzare quella che si può descrivere come una riscrittura oulipiana della Matematica».[7] Tra i membri fondatori figurano matematici del calibro di Henri Cartan, André Weil e Jean Dieudonné, animati dall'intento di rifondare l'intera matematica. Il progetto iniziale infatti è quello di redigere un trattato di Analisi Matematica per sostituire nell'insegnamento superiore i manuali esistenti manifestamente insufficienti ma, ben presto, questo obiettivo si è rivelato decisamente riduttivo. Così nel 1939 Bourbaki pubblica il trattato *Éléments de mathématique* (testo noto anche come *Traité*) concepito, come già denuncia l'uso del singolare nel titolo piuttosto che l'usuale *mathématiques*, secondo una visione unitaria e assiomatica per evitare una frammentazione della matematica in branche non comunicanti.

Le Lionnais considera la concezione bourbakista di larghe vedute; infatti «[Bourbaki] non ammette che una teoria sia chiusa su se stessa; non intende cristallizzare la matematica né fare un'opera definitiva; non pretende in alcun modo di coltivare la generalizzazione per se stessa, ma solo nella misura in cui risulti suscettibile di dare frutti»[8]. Caposaldo della matematica bourbakista è il *metodo assiomatico*, articolato sullo schema assioma-definizione-teorema, come si può leggere nella prima pagina degli *Éléments*.[9]

---

[5] «Depuis que j'ai eu le plaisir de faire votre connaissance, écrit Le Lionnais le 10 décembre 1943 […] J'ai, naturellement, réfléchi à notre conversation. Très séduit par l'une et l'autre des deux idées que vous m'avez proposées: 1) Les Mathématiques comme terme de la classification des Sciences, la Biologie en étant le début. 2) Application des Mathématiques à la description des faits historiques, et éventuellement à leurs prévision» [41, p. 120].
[6] Movimento artistico-letterario fondato a Parigi negli anni '20 del Novecento da André Breton.
[7] «Bourbaki was a group of French mathematicians who decided to perform what one may describe as an oulipian rewriting of Mathematics» [34, p. 126].

Le Lionnais e Queneau, entrambi letterati *e* matematici[10], sono interessati allo stato della Matematica contemporanea e dunque fortemente attratti dall'opera di Bourbaki che, secondo un'ipotesi di Roubaud, «quando l'Oulipo fu concepito, […] fornì un contro-modello al gruppo surrealista».[11]

Il piano programmatico di Bourbaki di rifondare la Matematica a partire dalla *Teoria degli insiemi* e ricorrendo solo al *metodo assiomatico* viene quindi esplicitamente importato nella fondazione dell'Oulipo da Le Lionnais il quale afferma senza mezzi termini: «Nel mio lavoro, Dada è stato all'origine del Terzo settore e Bourbaki dell'Oulipo»[12].

Per usare le parole dell'*oulipien* Olivier Salon: «Le Lionnais era un uomo di idee, un uomo di strutture, un uomo di tassonomia, molto più che uno scrittore. Per lui Bourbaki aveva svolto un compito importante di classificazione della matematica e di completa ricostruzione. Questo modello non poteva che imporsi in lui quando si è trattato di lavorare su *la contrainte*, nonostante le esitazioni del gruppo, che sono molto evidenti nei primi tre anni (e anche di più, diciamo fino all'arrivo di Roubaud nel 1966 e poi di Perec l'anno dopo), come si può vedere nei Resoconti di Jacques Bens sugli incontri dell'Oulipo nei primi anni [9]»[13].

Nell'intervista *Un Certain Disparate* con Jean-Marc Levy-Leblond e Jean-Baptiste Grasset [43] Le Lionnais racconta di aver conosciuto le idee bourbakiste mediante Enrique Freymann, direttore delle edizioni Hermann, casa editrice di Bourbaki, e dichiara di essere diventato iperbourbakista solo dopo aver superato una certa reticenza iniziale. Egli tuttavia rimprovera ai bourbakisti una certa aristocrazia e il chiaro disprezzo per l'insegnamento elementare e, cosa più importante, per la divulgazione, mentre si trova in linea con la loro scuola di pensiero per quanto riguarda l'avversione per le applicazioni della matematica. Citando Le Lionnais: «mi sono quindi convertito a Bourbaki e volevo, per il mio libro (*Les grands courants de la pensée mathématique*), avere un articolo di Bourbaki su Bourbaki»[14]. Le Lionnais afferma di esser stato messo in contatto con alcuni bourbakisti, tra cui Dieudonné, da Charles Ehresmann, uno dei padri fondatori della scuola bourbakista, e di aver faticato parecchio a convincerli a scrivere un articolo di diffusione scientifica. Tale contributo, dal titolo *L'architecture des Mathématiques* tradotto in inglese e pubblicato nel 1950 nell'*American mathematical monthly* [14][15], è il primo e l'unico lavoro non matematico sul progetto del gruppo[16] che ha sicuramente contribuito ad estendere oltreoceano la fama internazionale di Bourbaki. Le Lionnais riferisce di aver discusso con i bourbakisti a proposito della definizione che essi danno della matematica, ovvero se «la matematica è lo studio delle strutture che operano necessariamente sugli insiemi oppure se, facendo a meno degli insiemi, la matematica è lo studio delle strutture».[17] Le Lionnais osserva che «lo stile di pensiero di Bourbaki si colloca così nel grande sforzo moderno della formalizzazione della matematica. Impone ai suoi seguaci una volontà di purificazione quasi ascetica. Costringendosi a ricostruire la matematica solo passo dopo passo, si assicura una straordinaria solidità, ma evidentemente si condanna a perdere i benefici del rischio. Si riconosce l'impostazione caratteristica dei logici; non è sempre quella degli inventori»[18].

---

[8] «Bourbaki […] n'admet pas qu'une théorie soit fermée sur elle-même; qu'il n'entend pas figer les mathématiques ni faire une oeuvre définitive; qu'il ne pretend en aucune manière cultiver la généralisation pour elle-même, mais seulement dans la mesure où elle se révèle susceptible de porter des fruits» [18, p. 10].

[9] «Dai greci, chi dice matematica dice dimostrazione. Alcuni dubitano che al di fuori delle matematiche esistano dimostrazioni nel senso preciso e rigoroso che questo termine ha ricevuto dai greci e che si intende dare in questa opera. Si ha il diritto di dire che il significato del termine dimostrazione non è variato, poiché ciò che è stato una dimostrazione per Euclide, lo è tuttora ai nostri occhi; ed in epoche nelle quali tale nozione ha rischiato di perdersi e la matematica si è trovata in pericolo, è presso i greci che si è ricercato il modello. Ma a questa venerabile eredità si sono aggiunte, da un secolo, importanti scoperte. In effetti l'analisi del meccanismo di dimostrazione nei migliori testi di matematica ha permesso di liberare la struttura dal doppio punto di vista del vocabolario e della sintassi. Si arriva quindi alla conclusione che un testo di matematica sufficientemente esplicito può essere espresso in un linguaggio convenzionale comprendente solamente un piccolo numero di termini invariabili assemblati mediante una sintassi che consisterà in un piccolo numero di regole inviolabili. Un testo così concepito si dice formalizzato».

[10] Giusto per citare delle pubblicazioni di carattere prettamente matematico ricordiamo per Le Lionnais [17] e [20] e per Queneau [33] e le sue speculazioni sui *numeri iperprimi* di cui Le Lionnais riferisce in [22] e [23].

[11] «When the Oulipo was conceived, Bourbaki provided a counter-model to the Surrealist group» [34, p. 127].

[12] «Dans mon propre travail, Dada fut à l'origine du Troisième secteur, et Bourbaki de l'Oulipo» da una intervista di Ralph Rumney a Le Lionnais apparsa inizialmente in inglese nel 1975 su *Studio International* e poi ritradotta in francese da Patrice Cotensin, e pubblicata nel 1997 su *L'Échoppe* col titolo «Marcel Duchamp joueur d'échecs et un ou deux sujets d'y rapportant».

[13] «Le Lionnais était un homme d'idées, un homme de structures, un homme de taxonomie, bien plus qu'un écrivain. Pour lui, Bourbaki avait accompli une tâche majeure de classification des mathématiques, et de reconstruction complète. Ce modèle n'a pu que s'imposer à lui quand il s'est agi de travailler sur la contrainte, en dépit des hésitations du groupe, qui sont très claires dans les trois premières années (et même plus, disons jusqu'à l'arrivée de Roubaud, puis de Perec, en 1966 et 67), comme on le voit dans les Comptes rendus par Jacques Bens des réunions de l'Oulipo des premières années [9]» O. Salon, corrispondenza epistolare privata con le autrici (2020).

[14] «Je suis donc converti à Bourbaki et je voulais, pour mon livre (Les grands courants de la pensée mathématique), avoir un article de Bourbaki sur Bourbaki».

[15] In [43] Michèle Audin ipotizza che l'articolo di Bourbaki pubblicato ne *Les grands courants* sia stato scritto dal solo Dieudonné.

La genesi dell'articolo a nome di Bourbaki ne *Les grands courants* suggerisce che tra i due gruppi, oltre a intersezioni teoriche e metodologiche, vi furono anche contatti tra persone.

Le Lionnais racconta in [43] di avere contattato diversi bourbakisti. In particolare, oltre ai già citati Cartan, Dieudonné, Weil, Ehresmann e Godement, anche Jean Delsarte e Claude Chevallay. Queste interazioni si concretizzarono nei seguenti articoli per *Les grands courants*: *L'avenir des mathématiques* di Weil[19], *Les méthodes modernes et l'avenir des Mathématiques concrètes* di Godement e *David Hilbert (1862-1943)* e *Les méthodes axiomatiques modernes et les fondements des mathématiques*[20] di Dieudonné. Cruciale fu anche il ruolo svolto da Élie ed Henri Cartan, padre e figlio, coinvolti come intermediari tra Le Lionnais e Weil che risiedeva a San Paolo del Brasile tra il 1945 e il 1947.

Non solo Le Lionnais conosceva alcuni membri di Bourbaki di cui era coetaneo, ma anche Queneau fu per esempio ospite del Congresso di Amboise nel 1962 trattamento d'eccezione usualmente riservato solo alle reclute del gruppo. Queneau, estimatore delle opere di Bourbaki, in [30][21] dedica al *Traité* un'analisi concreta e attenta, intessuta della tipica ironia oulipiana come quando osserva che: «da lettura degli *Eléments*, dice l'opuscolo informativo, "non presuppone (…) in linea di principio, alcuna conoscenza matematica particolare, ma solo una certa abitudine al ragionamento matematico e un certo potere di astrazione". Naturalmente, questa frase non deve essere presa troppo alla lettera; e il "certo" può essere qualificato come una litote».[22] Secondo Queneau inoltre Bourbaki si rivolge «principalmente ai giovani che hanno una buona formazione matematica e che hanno la mente "aperta". […] D'altronde non bisogna credere che, meraviglia genetica, Bourbaki si sia generato da solo. Suo padre stesso, lo si può nominare, è Hilbert, e egli ebbe come nutrice, tra l'altro, Emmy Noether, e i suoi padrini furono quasi tutti degli stranieri».[23] In queste parole vi è il chiaro riferimento al periodo di stasi della matematica francese del primo Novecento, dovuta, secondo lui, sia alla prima guerra mondiale che all'influenza filosofica di Henri Poincaré. Pertanto «intorno agli anni 30, alcuni giovani matematici hanno preso coscienza del ritardo accumulato dalla matematica francese non solo nell'insegnamento […] ma anche nella ricerca; così hanno creato Bourbaki; e dopo il 1945, tre Francesi hanno ricevuto la medaglia Fields, l'equivalente del premio Nobel della matematica, Schwartz, Serre e Thom».[24]

---

[16] «François Le Lionnais had managed to persuade Bourbaki to write its first and only non-mathematical paper stating the project of the group» [34, p. 127].
[17] «[J'ai eu d'ailleurs des discussions avec eux à propos de cet article, et je me souviens leur avoir posé un problème sur la définition qu'ils donnaient de la mathématique, qu'ils ont discuté sans pouvoir le résoudre: ou bien,] la mathématique est l'étude d'ensembles munis de structures, ou bien la mathématique est l'étude des structures» [43] https://blogs.oulipo.net/Le Lionnais/2010/10/10/36-bourbaki/ (Accessed 22 April 2021).
[18] «Le style de pensée de Bourbaki se situe ainsi dans le grand effort moderne de formalisation des mathématiques. Il impose à ses partisans une volonté de purification presque ascétique. En s'obligeant à ne reconstruire les mathématiques que pas à pas, il s'assure une extraordinaire solidité, mais il se condamne évidemment à perdre les bénéfices du risque. On reconnaît la démarche caractéristique des logiciens; elle n'est pas toujours celle des inventeurs» [18, p. 9].
[19] A giudizio di Le Lionnais [43] l'articolo più difficile dal punto di vista matematico di tutto il volume.
[20] Questo articolo di Dieudonné non figura in [17] ma solo nell'edizione ampliata del 1962 edita da Librairie Scientifique et Technique Albert Blanchard, Paris.
[21] Tale contributo è successivamente ri-pubblicato come primo capitolo del suo libro *Bords. Mathématiciens, précurseurs, encyclopédistes* del 1963 edito da Hermann.
[22] «La lecture des Eléments, dit le prospectus, "ne suppose (…) en principe, aucune connaissance mathématique particulière, mais seulement une certaine habitude du raisonnement mathématique et un certain pouvoir d'abstraction". Naturellement, il ne faut pas prendre cette phrase trop au pied de la lettre; et le "certain" peut être qualifié de litote» [30, p. 6].
[23] «principalement aux jeunes qui on tune bonne formation mathématique et qui ont l'esprit "ouvert". […] Il ne faut pas croire d'ailleurs que, merveille génétique, Bourbaki se soit engendré lui-même. Son père même, on peut le nommer, c'est Hilbert, et il eut comme nourrice, entre autres, Emmy Noether, et ses parrains furent quasiment tous des étrangers» [30, p. 6].

Infine Queneau non perde mai di vista il legame tra matematica e arte sostenendo che: «utilità, bellezza, […] sono le due caratteristiche della matematica, quelle che la avvicinano all'arte e la differenziano da essa. Una teoria matematica vivente […] è sia bella che utile. E questo senza che ci sia contraddizione tra questi due aspetti».[25] Anzi, secondo Queneau, poiché non vi è contrapposizione tra l'ispirazione, l'intuizione e l'imposizione di *contrainte* [16, p. 343], il lavoro dell'Oulipo si sviluppa anche nella direzione di ridurre la distanza [29, p. 26], se non addirittura l'ostilità, tra letterati e matematici in vista della realizzazione di una sintesi perfetta tra scienza e umanesimo.

## Jacques Roubaud: un trait d'union?

Di primo acchito potrebbe sembrare che Roubaud, definito come il più bourbakista di tutti gli *oulipiens*, costituisca un vero e proprio trait-d'union tra i due *milieu*. In realtà per comprendere fino in fondo la sua complessa personalità è necessario non decontestualizzarla dalla sua vita[26]. Nasce il 5 Dicembre 1932 a Caluire-et-Cuire, vicino Lione, in una famiglia d'intellettuali: il padre Lucien, ex-normalista, è un professore di filosofia, la madre Suzanne Molino, una delle prime tre donne ad essere ammesse nel 1927 all'École Normale Supérieure, è una professoressa d'inglese e il fratello maggiore sarà un docente universitario di matematica. Come afferma Véronique Montémont «diventare un universitario è quasi un'abitudine familiare: questo contesto chiarisce che Roubaud abbia potuto trovare una via soddisfacente malgrado i suoi dubbi e permettersi di non scegliere tra le due sue passioni».[27] Il suo primo amore è certamente la poesia: i primi componimenti *Poésies juvéniles* risalgono al periodo 1942-44. Successivamente grazie ad un insegnante del liceo scopre i surrealisti: Éluard, Breton e Aragon. Conseguito il diploma, dopo vari tentennamenti, nel 1954-1955 decide di studiare matematica alla Sorbona, scelta che spiazza i suoi genitori. Infatti, se è vero che il suo desiderio di dedicarsi alla poesia è un punto fermo, è altrettanto vera la consapevolezza che tale attività non gli darà da vivere e, secondo l'autrice di *Jacques Roubaud: l'amour du nombre*, la matematica gli «consente nello stesso tempo una sorta di compartimento stagno delle due vite intellettuali, che lascia la pratica della poesia al di fuori dell'approccio universitario»[28]. Così Roubaud, studente sia di simpatizzanti di Bourbaki, come Gustave Choquet, che di collaboratori, come Laurent Schwartz, totalmente affascinato da quello che lui stesso definisce «una sorta di surrealismo matematico» [36, p. 393], diviene successivamente allievo del bourbakista Chevalley. Quest'ultimo, come riferisce Michèle Chouchan, durante il periodo di preparazione della tesi cerca di spronare il suo studente a scrivere una sorta di bozza di un capitolo sulle categorie, lavoro che però non è mai venuto al mondo.[29] Roubaud quindi non è un "militante di Bourbaki", piuttosto si avvale del bourbakismo per perseguire con determinazione l'essere poeta, come dichiara in *Poésie* (2000) «utilizzo il bourbakismo unicamente come manuale sul modo di comportarsi davanti ai problemi che pone la composizione della poesia»[30]. Contemporaneamente afferma anche che «l'idea della poesia come arte, come artigiano e come passione, come gioco, come ironia, come ricerca, come conoscenza, come violenza, come attività autonoma, come forma di vita, idea che è stata quella di molti poeti (quelli che preferisco) nella tradizione europea e, di

---

[24] «Vers les années 30, quelques jeunes mathématiciens ont pris conscience du retard pris par les mathématiques françaises non seulement dans l'enseignement (il ne cesse de l'être que depuis peu) mais encore dans la recherché; et ils créèrent Bourbaki; et depuis 1945, trois Français ont eu la médaille Field, sort de prix Nobel des mathématiques, Schwartz, Serre et Thom» [30, p. 7].

[25] «Utilité, beauté, voilà bien les deux caractères de la mathématique, ceux qui la rapprochent de l'art et l'en différencient. Une théorie mathématique vivant (et vrai en soi - mais c'est une autre histoire) est à la fois belle et utile. Et ceci sans qu'il y ait contradiction entre ces deux aspects» [30, p.5].

[26] Per una biografia dettagliata di Roubaud si rimanda, per esempio, a [24].

[27] «Devenir un universitaire est presque une habitude familiale: ce contexte explique que Roubaud ait pu trouver une voie satisfainte malgré ses doutes et se permettre de ne pas choisir entre ses deux passions» [24].

[28] «permettre en même temps un cloisonnement étanche des deux vies intellectuelles, ce qui laisse la pratique de la poésie en dehors de toute approche universitaire» [24].

[29] «Au moment où il préparait sa thèse, Roubaud a été sollicité par Claude Chevalley pour contribuer à une esquisse d'un chapitre sur les catégories, qui finalement, n'a jamais vu le jour» [15, p.124].

[30] «J'utilise le bourbakisme uniquement comme manuel sur la manière de se comporter devant les problèmes que pose la composition de la poésie».

recente, quella anche di Raymond Queneau, io l'ho fatta mia, e ne vedo un primo esempio tra i trovatori».[31]

Chiaramente l'impronta bourbakista figura e persiste tanto nella sua produzione matematica quanto in quella letteraria soprattutto a livello stilistico [24]. Se è vero infatti che, come vedremo, a differenza dell'Oulipo, Bourbaki si è ormai quasi auto-estinto, è anche vero che il suo *stile* di scrittura improntato alla chiarezza espositiva e compositiva gli è sopravvissuto e ormai fa parte del patrimonio culturale di chiunque scriva articoli di matematica. Un ulteriore livello di influenza consiste nella manifesta intenzione di fondare il suo progetto letterario sulla puntuale traslazione di nozioni matematiche nel testo poetico [24]. Un titolo come ∈ (*Signe d'appartenance*) è un esplicito rimando alla Teoria degli insiemi[32]. Tale raccolta di poesie, concepita fin dal 1961 come reazione al primo avvenimento doloroso della sua vita, la morte violenta del fratello Jean-René, diviene il primo atto concreto della realizzazione di ciò che egli chiama il suo Progetto[33]. Per Roubaud il disincanto da Bourbaki è direttamente proporzionale alla fascinazione che aveva subito negli anni giovanili.[34] Ad un certo punto egli diventa critico sul *Traité*: «Ci sono inoltre un sacco di argomenti su cui Bourbaki è confuso o rimane silenzioso. Non ci si rende conto che è una presentazione spesso ingannevole. A dire il vero ho avuto in mano prima ciò che è divenuto il primo volume di topologia generale, che continuo a pensare sia il migliore di tutti. Almeno agli inizi di Bourbaki, la sua lingua era molto bella e infinitamente seducente. Ma poi, lo stile si è intasato di tecnicità».[35] Come ipotizza Montémont quasi certamente l'opinione di Weil, secondo cui «esistono solo i grandi matematici, i "cercatori", che sono dei geni precoci e sono all'apice del loro talento tra i venti e i trent'anni [e] coloro che non hanno questa fortuna sono paragonati ad una "cassa di risonanza" ad uso degli altri»[36], ha lasciato un segno in Roubaud. Considerato che egli ha iniziato tardivamente ad occuparsi di matematica, contribuendo in modo significativo ma non geniale alla ricerca matematica, molto probabilmente ciò ha influito nel suo proposito di mettere la sua esperienza scientifica al servizio della composizione poetica. In *Poésie, etcetera: ménage* [37], saggio della sua maturità letteraria, ritorna sui rapporti tra Oulipo e Bourbaki rivelando che «si può anche dire che l'Oulipo è un omaggio a Bourbaki, un'imitazione di Bourbaki. Allo stesso tempo è, in maniera meno evidente, una parodia di Bourbaki, se non una profanazione. Il progetto oulipiano, che "traduce" la visione e il metodo bourbakista nel campo delle arti del linguaggio è altrettanto serio, ambizioso, ma non settario, e non persuaso della validità della sua proposta in opposizione a qualsiasi altro approccio».[37]

---

[31] «L'idée de poésie comme art, comme artisanat et comme passion, comme jeu, comme ironie, comme recherche, comme savoir, comme violence, comme activité autonome, comme forme de vie, idée qui fut celles de bien de poètes (ceux que je préfère) dans la tradition européenne, et tout récemment encore celle de Raymond Queneau, je l'ai faite mienne, et j'en vois l'exemple premier chez les troubadours» [*La fleur inverse, essai sur l'art formel des troubadours*. Ramsey, Paris (1986) p. 17].

[32] In [24] l'autrice sottolinea come *Signe d'appartenance*, composta durante la stesura della tesi di matematica, sia «da raccolta in cui *la contrainte* sembra più forte, ma curiosamente è stata scritta quando Roubaud non si era ancora unito al gruppo né conosceva personalmente Queneau. Il giovane poeta ha incontrato quest'ultimo quando ha scelto di sottoporgli la sua raccolta per la pubblicazione presso Gallimard».

[33] Scrive Roubaud in *Le grand incendie de Londres*: «Je devais poursuivre ailleurs ma préparation au projet: dans la mathématique, dans la poésie, dans une grande sévérité d'existence. L'austérité parfois érémitique qui se montrait nécessaire était comme fonctionnellement imposée par une recherche simultanée de voies dans les deux directions, duales et antagonistes en apparence, de la mathématique et de la poésie».

[34] Egli riferisce che il contenuto degli *Éléments de géométrie algébrique* di Grothendieck e Dieudonné «coulait dans ma cervelle comme un miel, que dis-je? un nectar, une ambroisie intellectuelle. Je ne peux m'en souvenir sans stupéfaction» [*Mathématique*. Seuil, Paris (1997)].

[35] «Il y a d'ailleurs des tas de sujets sur lesquels Bourbaki est confus ou demeure silencieux. On ne se rend pas compte que c'est une présentation souvent fallacieuse. En fait j'ai d'abord eu entre les mains était le premier tome de topologie générale, que je persiste à penser être le meilleur de tous. Au moins dans les débuts de Bourbaki, sa langue était très belle, et infiniment séduisante. Mais ensuite, le style s'est engorgé de technicité» [15, p.124].

[36] «n'existent que les grands mathématiciens, les "trouveurs", qui sont des génies précoces et sont à l'apogée de leur talent entre vingt et trente ans. Ceux qui n'ont pas cette chance sont comparés à une "caisse de résonance" à l'usage des autres» [24].

[37] «On peut dire aussi que l'Oulipo est un hommage à Bourbaki, une imitation de Bourbaki. En même temps il est, de manière non moins évidente, une parodie de Bourbaki, sinon une profanation. Le projet oulipien, qui "traduit" la visée et la méthode bourbakiste dans le domaine des arts du langage est également sérieux, ambitieux, mais non sectaire, et non persuadé de la validité de sa démarche à l'exclusion de toute autre approche» [37, p. 201].

Introdotto da Queneau nel 1966 entra a far parte dell'*Ouvroir* e nel 1981 fonda l'*Atelier de Littérature Assistée par les Mathématiques et les Ordinateurs* (ALAMO). Nel 1983 un'altra vicenda dolorosa segna la vita del "poeta formalista"[38]: l'improvvisa morte della giovane moglie Alix Cléo, sposata qualche anno prima. Dopo un periodo di silenzio di oltre due anni e mezzo, Roubaud inizia a ideare un ciclo di prosa autobiografica (*traités de mémoire*), la cui prima pubblicazione nel 1989, *Le grand incendie de Londres*, descrive il fallimento del suo Progetto.

Grazie al fatto che l'Oulipo permette a Roubaud di attuare concretamente in una sola impresa artistica la sua duplice passione per la letteratura e la matematica, il suo ruolo attivo da *oulipien* costituisce un terreno fertile per la sua produzione umanistica: non solo opere di prosa e poesia – a nome singolo o meno, tra i ranghi dell'Oulipo o meno – ma anche saggi, lavori teatrali e di critica letteraria.

Oltre ai numerosi contributi nei volumi collettivi [27, 28] e ai circa 30 fascicoli della *Bibliothèque oulipienne* di cui è autore, volendo citare alcune opere ricordiamo per esempio il cosiddetto *Ciclo di Ortensia*. Si tratta di un ciclo di sei romanzi - di cui solo tre sono stati pubblicati[39] - in cui i rimandi allo schema della sestina lirica sono evidenti tanto nel piano dell'opera quanto nella trama. Inoltre poiché, come afferma Queneau «il sonetto è una delle […] cavie preferite»[40] degli *oulipiens*, Roubaud si è cimentato, da matematico, nella speculazione teorica sottesa alla composizione di quenine[41] e, da poeta, nella creazione delle stesse. Nel primo contesto a lui sono dovuti alcuni risultati[42], come il teorema che porta il suo nome, sulla caratterizzazione dei cosiddetti *numeri di Queneau* ovvero quei numeri $n$ per i quali è possibile comporre una *n-ina*. Sul fronte della sperimentazione poetica invece vale la pena ricordare la più recente raccolta di 36 *pharoïnes*, ossia quenine composte utilizzando una particolare permutazione che lo stesso Roubaud ha ideato, pubblicate nel volume *Les fastes* che accompagna un'esposizione dell'artista Jean-Paul Marcheschi.[43]

Intellettuale dai molteplici interessi, dalla produzione ricca e variegata, Roubaud dopo la morte di Le Lionnais e Queneau rappresenta uno dei pilastri del gruppo oulipiano e ha contribuito ad accrescerne la fama al pari di autori come Perec e Calvino.

## Lo strutturalismo bourbakista nel progetto dell'Oulipo

In [25] Montémont individua una comune radice tra Bourbaki e Oulipo nell'attitudine tipica dello Strutturalismo di concepire ogni oggetto di studio come una struttura scomponibile in elementi e unità, il cui valore funzionale è determinato dall'insieme dei rapporti fra ogni singolo livello dell'opera e tutti gli altri.

Per quanto riguarda l'Oulipo infatti, benché ne *Le second Manifeste* [28] Le Lionnais prenda inequivocabilmente le distanze dall'approccio strutturalista[44] affermando che l'Oulipo ha una nozione più sviluppata di cosa sia una struttura ed essa deriva proprio da Bourbaki, allo stesso tempo afferma che l'attività principale dell'Oulipo è investigare l'applicazione di alcune strutture matematiche ai diversi livelli della produzione letteraria grazie all'imposizione di vincoli e regole (*contrainte*)[45]. L'Oulipo si è inoltre dato il duplice compito di elencare tutti i vincoli esistenti e l'invenzione di nuovi vincoli. Nelle ricerche che l'Oulipo si propone di intraprendere infatti coesistono due tendenze principali, due *LiPos* [Lescure in [28] p. 33]. Una orientata all'Analisi, l'*Anoulipisme*, è votata alla scoperta ossia a rintracciare nello studio di opere del passato l'utilizzo di *contrainte* preesistenti alla loro teorizzazione da parte dell'Oulipo, scherzosamente detti *plagi per anticipazione*. L'altra orientata alla Sintesi, il *Sintoulipisme*, è votata all'invenzione di nuovi vincoli matematico-letterari.

---

[38] Tale appellativo, che Roubaud stesso evoca in [37], non è stato sempre utilizzato nella sua accezione positiva dato che le sue opere hanno spesso suscitato reazioni contrastanti ed egli è considerato un autore dal carattere controverso.

[39] *La belle Hortense* e *L'enlèvement d'Hortense* sono stati pubblicati da Ramsay, rispettivamente nel 1985 e nel 1987, mentre *L'exil d'Hortense* è stato pubblicato nel 1990 da Seghers.

[40] «de sonnet est un de […] cobayes préféerés» [29, p. 26].

[41] Per una definizione di *quenina* o *n-ina* si rimanda al paragrafo successivo.

[42] Roubaud J.: Un problème combinatoire posé par la poésie lyrique des troubadours. Mathématiques et Sciences Humaines 27, 5-12 (1969) e Roubaud J.: N-ine autrement dit quenine (encore). La Bibliothèque Oulipienne 66, (1993).

[43] Marcheschi, J.-P. and Roubaud J.: Les Fastes. Lienart et Musée départemental de préhistoire d'Île de France (2009).

[44] «La très grande majorité des œuvres OuLiPiennes qui ont vu le jour jusqu'ici se place dans une perspective syntaxique structurEliste (je prie le lecteur de ne pas confondre ce dernier vocable – imaginé à l'intention de ce Manifeste – avec StructurAliste, terme que plusieurs d'entre nous considèrent avec circonspection» [28, p. 19].

[45] «L'activité de l'Oulipo et la mission dont il se considère investi pose le(s) problème(s) de l'efficacité et de la viabilité des structures littéraires (et, plus généralement, artistiques) artificielles» [28, p. 20].

D'altro canto, nella concezione di Bourbaki la matematica appare come una sorta di contenitore di forme astratte, le strutture matematiche, e il metodo assiomatico stesso che, secondo Dieudonné consiste nello studio logico delle «*relazioni* tra enti e non la loro *natura*», opera mediante il concetto di strutture che sono gli strumenti del matematico.

È quindi evidente come l'anima e la vocazione strutturalista di Bourbaki abbia indubbiamente influenzato Le Lionnais, meno attraverso la sua matematica, quanto piuttosto attraverso l'enfasi posta sulle strutture e sulla loro definizione a partire da un assunto assiomatico. Così, secondo Roubaud in [34, p. 127], se Bourbaki lavora con le strutture, Oulipo lavora con i vincoli che sono l'equivalente oulipiano delle strutture bourbakiste e per quanto riguarda il metodo afferma: «potremmo dire che il metodo oulipiano *imita* il metodo assiomatico, che ne è una trasposizione, un trasporto nel campo della letteratura».[46] Da ciò segue che un testo scritto secondo i dettami della letteratura potenziale è l'equivalente di un teorema tanto più che il lavoro proprio dell'Oulipo non è produrre testi letterari ma, come già detto, inventare, scoprire o riscoprire *contrainte* che fungano da strumento creativo per gli scrittori: «La *LiPo* si propone di definire delle *strutture*. All'interno di queste strutture, una personalità surrealista, o classica, o romantica, o tutto ciò che si vuole, può esprimersi e quindi produrre testi surrealisti, classici, romantici, etc. Ma la *LiPo* in sé non è né questo né quello».[47] Per Le Lionnais «il vincolo si esaurisce nel gesto della sua scoperta, nella sua definizione. In definitiva, non è necessario che nemmeno un singolo testo venga scritto soddisfacendo un vincolo precedentemente formulato».[48] Per questo egli non ha mai posto l'enfasi sulla creazione di testi ma solo sulla corretta definizione dei vincoli ai fini dell'uso sistematico di strutture per la creazione letteraria. Jacques Duchateau ricorda che l'Oulipo «dovrebbe avere il progetto di esplorare la nozione di struttura in letteratura, indipendentemente dall'uso che se ne può fare».[49]

In breve la *LiPo* «vuole essere, come Bourbaki in matematica, il supporto di una rivoluzione nella concezione stessa che si ha dell'atto di scrivere».[50]

In questo senso, sebbene Le Lionnais abbia composto con un linguaggio molto formale e matematico tutti i manifesti oulipiani, con il terzo sembra volersi spingere oltre delineando per il gruppo un progetto più marcatamente matematico ovvero classificare tutte le possibili strutture della *LiPo*: «Un programma di costruzione di tutte le strutture letterarie possibili – grande vocazione dell'OuLiPo – si articola in tre fasi che richiedono solo l'esercizio di un certo talento durante un tempo inferiore all'eternità. La prima fase pone un problema teorico, [apparentemente] semplice come la scoperta di tutti i possibili teoremi matematici o la preparazione di tutte le possibili sostanze chimiche. La seconda fase, anch'essa teorica, mira a estrarre da questo cumulo le strutture di riconosciuta efficacia. Essenzialmente pratica, la terza fase ambisce a portare alla Soglia dell'Opera».[51] Gli *oulipiens* avranno dunque il compito di compilare e completare *Le Grand Tableau*[52] ovvero una matrice a doppia entrata in cui ogni colonna corrisponde a una

---

[46] «On pourrait dire que la méthode oulipienne *mime* la méthode axiomatique, qu'elle en est une transposition, un transport dans le champ de la littérature» [36, p. 404].

[47] «La LiPo se propose de définir des *structures*. À l'intérieur de ces structures, une personnalité surréaliste, ou classique, ou romantique, ou tout ce que l'on veut, peut s'exprimer et donner ainsi des textes surréalistes, classiques, romantiques, etc. Mais la LiPo elle-même n'est ni ceci, ni cela» [9, p. 143].

[48] «Pour Le Lionnais, enfin, la contrainte s'épuisait dans les geste de sa découverte, dans sa définition. À la limite, il n'était pas nécessaire même qu'un seul texte soit écrit satisfaisant à une contrainte oulipiennement formulée» [35, p. 81]

[49] «[L'Oulipo] devait avoir le projet d'explorer la notion de structure dans la littérature, indépendamment de l'utilisation susceptible d'en être faitre» *Prefazione* a [9, p. 15].

[50] «Elle se veut, comme Bourbaki en mathématiques, le support d'une révolution dans la conception même que l'on se fait de l'acte d'écriture» [13, p. 255].

[51] «Un programme de costruction des toutes les structures littéraisres possibles - vocation majeure de l'OuLiPo - passe par trois phases qui ne requièrent pas que l'exercice d'un certain talent durant un temps inférieur à l'éternité. La première phase pose un problème théorique, aussi simple [En apparence] que la découverte de tous les théorèmes de mathématiques possibles ou que la fabrication de toutes les substances chimiques possibles. La seconde phase, également théorique, vise à extraire de cet amoncellement des structures d'une efficacité reconnue. Essentiellement pratique, la troisième phase se donne pour ambition de conduire au Seuil de l'Œuvre.» [26, pp.798-799].

[52] *Le Grand Tableau* descritta da Le Lionnais nel Terzo Manifesto corrisponde a *La Table de Queneleiev* progettata da Queneau nel 1973 e pubblicata in [27, p. 73]. Una versione ampliata e aggiornata è disponibile all'indirizzo web https://www.fatrazie.com/jeux-de-mots/recreamots/253-table-de-queneleiev (Accessed 22 April 2021). Una seconda

struttura matematica e ogni riga a un oggetto letterario in modo che ciascuna casella sia definita dall'azione di una struttura matematica su un oggetto letterario. Più concretamente la prima fase del programma riguarderà la determinazione delle colonne e delle righe e il riempimento delle caselle; nella seconda fase, grazie a criteri che secondo Le Lionnais sono ancora per lo più da definire, si dovranno associare più caselle e studiare le loro giustapposizioni. Egli sceglierà di riferirsi a queste ultime con il termine *armature*, emendando i termini *struttura* e *forma fissa*[53] e ricordando così le "strutture madri" proposte da Bourbaki.
Nel 1976, in [43], Le Lionnais confessa di essere in certa misura deluso da quanto sino a quel momento ha realizzato l'Oulipo e affettuosamente accusa i membri del gruppo, soprattutto quelli più giovani, di essere più attratti dal clima allegro e familiare delle riunioni che interessati a ideare strutture potenti, efficaci e feconde, cosa in cui si sono spesi a suo dire poco frequentemente e poco incisivamente. Egli tuttavia afferma che «se mi viene chiesto di sapere se, a partire da strutture artificiali a priori, possiamo ritrovarci con qualcosa di vivo, la mia risposta è di un moderato ottimismo - che mi caratterizza in molti campi - è un'ipotesi, una scommessa»[54]. A tal fine ha concepito il Terzo Manifesto con lo scopo di «rimettere in sella con forza l'Oulipo» e come un lascito intellettuale non solo «ai membri dell'Oulipo ma a persone che verranno tra venti o quarant'anni e che potranno realizzare qualcosa in questo campo».
Le Lionnais annuncia sinteticamente per la prima volta la preparazione di un terzo manifesto durante la riunione del 28 ottobre 1975; altre menzioni e annunci si rincorrono nel '76, nel '78 e nell'80 creando un'attesa febbrile intorno alla sua pubblicazione la quale avverrà solo postuma [26] con il titolo *Le troisième manifeste. Prolégomènes à toute littérature future*. Il documento risulta però giusto l'introduzione a quello che sarebbe potuto essere un Manifesto e tale incompiutezza, suggerisce Camille Bloomfield, da un lato contribuisce ad alimentare il mito dell'Oulipo dall'altro forse renderà un giorno «possibile scrivere veramente, collettivamente, questo Terzo Manifesto, ma è più probabile che nella sua forma attuale, misteriosa, eccitante quasi, questo spettro di testo sia il più potente stimolo per ossessionare gli oulipiani, condannati (o incoraggiati) a continuare la ricerca teorica che è rimasta in forma di bozza»[55].
Alle "rimostranze" espresse da Le Lionnais sembra fare eco Roubaud che, in *La mathématique dans la méthode de Raymond Queneau*[56] [36]e in *Deux principes parfois respectés par les travaux oulipiens* [27, p. 90], propone quella che Natalie Berkman in [12] definisce «la seconda teorizzazione matematica del lavoro dell'Oulipo». In effetti in [36], dopo un'iniziale analisi degli aspetti matematici, dei metodi e dei risultati propri dei lavori

---

matrice nota come *La Tollé* (*Table des Opérations Linguistiques et Littéraires Élémentaires*) è stata concepita da Bénabou [8, pp. 101-106] il quale, nelle intestazioni delle colonne, ha sostituito le "strutture matematiche" della tavola di Queneau (*longuer*, *nombre*, *ordre* e *nature*) con "operazioni matematiche" (*déplacement et forme*, *substitution*, *addition*, *soustraction*, *multiplication (répétition)*, *division*, *prélèvement* e *contraction*). A Bénabou è dovuta anche *Une liste de contraintes oulipiennes* in cui in modo non esaustivo sono presentati alcuni dei *vincoli* oulipiani di cui, in certi casi, è data una definizione e se possibile un esempio (https://www.oulipo.net/fr/une-liste-de-contraintes-oulipiennes, Accessed 22 April 2021).

[53] «Toutes les potentialités que l'Oulipo […] a l'ambition d'explorer sont contenues dans ce que nous appellerons le GRAND TABLEAU. Aux imperfections près il constituera la base des deux premières phases et apportera d'heureuses suggestions aux volontaires de la troisième phase. Ce GRAND TABLEAU sera un quadrillage à double entrée, chaque colonne correspondant à une *structure mathématique*, chaque rangée à un *objet littéraire*. Chaque case sera donc définie par l'action d'une structure mathématique sur un objet littéraire. La détermination des colonnes et des rangées et le remplissage des cases constituent la première phase. La seconde phase consistera à essayer (en vertu de critères dont plusieurs restent à découvrir) d'associer plusieurs cases et à voir ce que donnent ces juxtapositions que nous appellerons des ARMATURES […]. Performance comparable à celle de la Nature (prolongée depuis peu par les chimistes) qui - à partir des molécules (les cases du GRAND TABLEAU) - fabrique tous ces mélanges que sont les substances naturelles (et artificielles). Les préoccupations de la troisième phase ne sont pas sans analogies avec celles qui hantent les chefs d'entreprises: conquête de certains marchés, recours à des technologies de pointe, rentabilité» [26, p. 799-800].

[54] «A cause de cela, nous n'avons pas abouti à des structures importantes et le but de mon troisième manifeste est de remettre l'OULIPO vigoureusement en selle. Dans ce troisième manifeste je ne m'adresse plus aux membres de l'OULIPO mais à des gens qui viendront dans vingt ou quarante ans et qui pourront réaliser quelque chose dans ce domaine. Je suis convaincu que c'est possible. Si on me pose la question de savoir si, de structures à priori artificielles on peut aboutir à quelque chose de vivant, ma réponse est d'un optimisme modéré – qui me caractérise dans pas mal de domaines – c'est une hypothèse, un pari». [43, punto 16] https://blogs.oulipo.net/Le_Lionnais/2010/10/10/16-le-troisieme-manifeste/ (Accessed 22 March 2021).

[55] «Peut-être les pistes ouvertes […] permettront-elles un jour d'écrire réellement, collectivement, ce Troisième Manifeste, mais il est plus probable que sous sa forme actual, mystérieuse, excitante presque, ce spectre de texte n'en soit que plus stimulant pour hanter les oulipiens, condamnés (ou incoraggagés) à poursuivre des recherches théoriques restées à l'état d'ébauche» [13, p. 366].

[56] Inizialmente pubblicato sulla rivista *Critique* nel 1977 è poi incluso nel 1982 nell'*Atlas* [27].

di Queneau, «si arriva alla creazione [...] dell'Oulipo, in cui la strategia precedentemente descritta diventa esplicita, sistematica e collettiva».[57] Il saggio ricalca, nel linguaggio e nella forma, la struttura degli articoli scritti dai bourbakisti grazie a un susseguirsi di assiomi, congetture, enunciati e proposizioni. Tre proposizioni, ad esempio, riassumono le caratteristiche del lavoro oulipiano che è naïf, divertente e artigianale[58]. Dopo i paragrafi intitolati *La contrainte*, *L'anti-hasard* e *La méthode axiomatique* in cui, ancora una volta, viene *enunciato*[59] lo stretto legame tra il metodo di Bourbaki e quello dell'Oulipo, il paragrafo *Les "Fondements de la Littérature"* riprende l'omonimo lavoro di Queneau *Les fondements de la littérature d'après David Hilbert* in cui l'autore, ispirandosi al celebre trattato del 1899 *Grundlagen der geometrie* di David Hilbert, presenta «[un approccio assiomatico alla] letteratura sostituendo nelle proposizioni di Hilbert le parole "punti", "linee", "piani", rispettivamente con: "parole", "frasi", "paragrafi"».[60] Infine, è la volta del paragrafo *Les structures* il cui incipit marca la distanza dallo Strutturalismo e la vicinanza a Bourbaki: «La struttura nel suo significato queniano e oulipiano non ha che una relazione minimale con lo "strutturalismo". Idealmente (come il vincolo rispetto all'assioma), essa si riferisce alla struttura bourbakista: l'oggetto, nel caso matematico, è un insieme con qualcosa "dentro" (delle leggi in algebra; degli intorni in topologia...); nel caso oulipiano l'oggetto è linguistico e la sua struttura è un modo di organizzarlo. Questa struttura soddisferà a una o più condizioni: assiomi in un caso, vincoli nell'altro; così un insieme munito di una legge di composizione avrà una struttura di monoide se questa legge soddisfa all'assioma di associatività; un testo avrà una struttura lipogrammatica se soddisfa al vincolo avente lo stesso nome».[61]

La matematica bourbakista chiaramente non è l'unica alla quale gli *oulipiens* hanno fatto riferimento nelle loro sperimentazioni al crocevia fra matematica e letteratura.

La letteratura potenziale, come dichiara lo stesso Le Lionnais in [23], si inserisce nell'ambito di un utilizzo della matematica che risale alla poesia provenzale medioevale[62]. Un esempio di 'struttura matematica letteraria' che ha catalizzato l'interesse dei membri dell'Oulipo è infatti quella della *sestina lirica* introdotta dal trovatore occitano Arnaut Daniel nel XII secolo [3]. Essa è caratterizzata dal fatto che ogni parola della prima strofa, formata da sei versi, ritorna nelle strofe successive secondo una fissata permutazione: da ciò l'idea di poter comporre delle *n*-ine con *n* diverso da sei - anche dette *quenine* in onore di Queneau dietro proposta di Roubaud - utilizzando i metodi propri della teoria dei gruppi [38]. Naturalmente le ricerche sulla generalizzazione della sestina non sono l'unico esempio di matematizzazione delle Lettere in cui l'Oulipo si è speso. Scorrendo le pagine delle opere collettive [26], [27] e [28] è facile ottenere una panoramica sui temi trattati e sui tentativi realizzati tra cui palindromi, lipogrammi, il metodo $M\pm7$, l'analisi matriciale del linguaggio, il meccano, dalla sestina alla quenina, la relazione '$x$ prende $y$ per $z$', le strutture, la combinatoria e così via. Quest'ultima prediletta dagli *oulipiens* non rientra invece fra gli interessi dei bourbakisti.

Dunque, nonostante ci siano circa 25 anni di distanza tra la nascita di Bourbaki (1935) e quella di Oulipo (1960), «l'intersezione epistemologica tra i due gruppi è reale» [25]. Chiaramente l'influenza non è reciproca tanto più che, nelle cognizioni attuali, nessun bourbakista pare si sia mai avvicinato seriamente all'Oulipo.

---

[57] «On en arrive alors à la création [...] de l'Oulipo, où la stratégie précédemment décrite devient explicite, systématique, et collective» [36, p. 399].

[58] «Proposition 8: Le travail oulipien est naïf. [...] Proposition 9: Le travail oulipien est amusant. [...] Proposition 10: Le travail oulipien est artisanal» [36, p. 399-400].

[59] «Proposition 14: Une contrainte est un axiome d'un texte. Proposition 15: L'écriture sous contrainte oulipienne est l'équivalent littéraire de l'écriture d'un texte mathématique formalisable selon la méthode axiomatique» [36, p. 404].

[60] «M'inspirant de cet illustre exemple, je présente ici une axiomatique de la littérature en remplaçant dans les propositions d'Hilbert les mots "points", "droites", "plans", respectivement par: "mots", "phrases", "paragraphes"» [32].

[61] «La structure dans son acception quenienne et oulipienne n'a qu'un rapport minimal avec les "structuralisme". Idéalement (comme la contrainte par rapport à l'axiome), elle se réfère à la structure bourbakiste: lìobject, dans le cas mathématique, est un (des) ensembles (s) avec quelque chose "dessus"(des lois en algèbre; des voisinages en topologie...); dans le cas oulipien l'object est linguistique et sa structure est un mode d'organisation. Cette structure satisfera à une ou plusieurs conditions: axiomes dans un cas, contrainte dans l'autre; ainsi un ensemble muni d'une loi de composition aura une structure de monoide si cette loi satisfait à l'axiome d'associativité; un texte aura une structure lipogrammatique s'il satisfait à la contrainte du même nom» [36, p. 410].

[62] «Je me décidai à proposer à Raymond de créer un atelier ou un séminaire de littérature expérimentale abordant de manière scientifique ce que n'avaient fait que pressentir les troubadours, les rhétoriqueurs, Raymond Roussel, les formalistes russes et quelques autres» in [23, p. 77].

## Analisi comparativa di due realtà parallele

Roubaud ha il merito di aver elencato alcune delle caratteristiche dell'Oulipo che a suo avviso sono in certa misura riconducibili all'influenza di Bourbaki [34, p. 128].

*Atto di nascita: luogo e nome.* Entrambi i gruppi hanno giocato molto sui propri aspetti leggendari e quasi mitologici. Tuttavia se da un lato Bourbaki, pur muovendosi nel solco della cultura normalista promotrice di scherzi e facezie, a causa di essa ha forse finito col sottolineare e affermare la propria natura elitaria ed esclusiva, in Oulipo la vena giocosa e burlesca – ereditata anche dal Collegio di Patafisica[63] – è prevalsa a beneficio dei singoli e della collettività tanto all'interno quanto verso l'esterno. Come dire che, pur coltivando entrambi il rigore compositivo ed espositivo delle loro opere[64], i bourbakisti si "prendevano molto sul serio" mentre gli *oulipiens* molto poco!
In questo ambito ricadono sicuramente la scelta del luogo in cui tenere la prima riunione fondazionale e la decisione del nome del gruppo.
L'idea di una possibile sperimentazione matematica nel processo creativo letterario è concepita da Le Lionnais nel corso dei suoi studi universitari. Tale idea, dopo averlo continuamente affascinato per molti anni, poté realizzarsi grazie alle conversazioni con Queneau, in particolare durante quelle svoltesi a Cerisy-la-Salle nel settembre 1960[65]. Giovedì 24 novembre 1960, nella cantina del ristorante *Au Vrai Gascon* di Parigi, Le Lionnais e Queneau fondano il *Seminaire de Littérature Expérimentale* (Selitex), di fatto una Sottocommissione del Collegio di Patafisica. Meno di un mese dopo, il 19 dicembre il Selitex viene ribattezzato «fino a nuovo ordine» [9, p. 28] *Ouvroir de littérature potentielle* (OLiPo) e infine, nel corso della riunione del 13 gennaio 1961, definitivamente consacrato come OuLiPo[66] grazie a una brillante intuizione di Albert-Marie Schmidt (1901-1966), uno dei dieci membri fondatori insieme a Noël Arnaud (1919-2003), Jacques Bens (1931-2001), Claude Berge (1926-2002), Jacques Duchateau (1929-2017), Latis *alias* Emmanuel Peillet (1913-1973), Jean Lescure (1912-2005) e Jean Queval (1913-1990).
Uno scenario simile vede la nascita di Bourbaki: il 10 dicembre 1934 intorno a un tavolo del *Capoulade* - rinomato café del *V arrondissement* parigino il cui menu proponeva piatti tipici della cucina alverniate[67] – si riuniscono per la prima volta Henri Cartan, René de Possel, André Weil, Jean Dieudonné, Claude Chevalley, Jean Coulomb e Jean Delsarte.
La questione onomastica è invece un po' più complessa. Pare che la scelta del nome dato al gruppo, nel segno di quell'atteggiamento burlesco di cui si è già detto, sia avvenuta nel corso del primo Congresso tenutosi a Besse-en-Chandesse nel luglio del 1935 e sia riconducibile al cognome del generale franco-cretese Charles Denis Bourbaki il quale combatté senza successo la guerra franco-prussiana (1870-1871). Poiché i membri del gruppo usavano chiamarsi tra loro per cognome [5], il problema del nome si pose solo quando, nell'autunno del 1935, fu necessario redigere una nota biografica[68] per accompagnare la pubblicazione di un articolo sui *Comptes Redus de l'Académie des Sciences*. Éveline Gillet (moglie di de Possel e poi di Weil) propose Nicolas poiché a suo dire suonava bene col cognome russo di Bourbaki [7]. Alcuni anni dopo in una nota biografica[69], presumibilmente del 1949, le origini del fittizio matematico si fecero risalire alla stirpe dell'illustre generale.

---

[63] Il Collegio di Patafisica era nato a Parigi l'11 maggio 1948.
[64] In [19] Le Lionnais sottolinea come l'Oulipo sia nato «nel duplice segno di una grande fantasia e di un non meno grande rigore» [42, *Conclusion*].
[65] Presso il Centre Culturel International de Cerisy-la-Salle si svolgeva una *décade* di colloqui dedicata a "Raymond Queneau ou une nouvelle défense et illustration de la langue française" organizzata da Georges-Emmanuel Clancier e Jean Lescure [6, *Préface*].
[66] «Il Sig. Latis fa subito notare che l'abbreviazione *OuLiPo* è preferibile a *OLiPo*. L'osservazione sembra densa di significato e la proposta che essa implica viene adottata» («M. Latis fait aussitôt remarquer que l'abréviation *OuLiPo* est preferable à *OLiPo* La remarque paraît pleine de sens et la propositionqu'elle sous-entend est adoptée») [9, p. 33].
[67] Michèle Audin racconta della nascita di Bourbaki e del suo legame con la regione francese dell'Alvernia in [5].
[68] Tale nota biografica fu letta agli accademici da Élie Cartan nella seduta del 23 dicembre 1935.

La storica Liliane Beaulieu ha sottolineato come la scelta di un *café* non sia certo sorprendente per l'epoca. Per entrambi i gruppi, incontrarsi in questo tipo di luogo sembra una dichiarazione programmatica emblematica della volontà di collocarsi al di fuori di realtà accademiche o letterarie o editoriali alle quali spesso sono affiliati i singoli membri ma da cui invece il gruppo nella sua interezza desidera prendere le distanze [6].

*Segretezza*. Bourbaki alimenta l'alone di mistero intorno alle sue origini e alle sue attività non solo celandosi dietro al nome fittizio, ma anche imponendo a ogni membro l'assoluto riserbo al punto che nessun bourbakista può ammettere pubblicamente di esserlo. L'Oulipo nasce come «una specie di società segreta»[70] ma dopo un primo decennio "sotto traccia" diviene noto al grande pubblico da un lato grazie ai notevoli successi editoriali di opere come *La Disparition* e *La Vie mode d'emploi* di Perec o *Zazie dans le métro* di Queneau, dall'altro in tempi più recenti per attività come i celebri "Jeudis de l'Oulipo", sessioni pubbliche di presentazione dei vari lavori del gruppo che hanno luogo un giovedì al mese[71] presso l'Auditorium della *Bibliothèque Nationale de France* a Parigi. Al passo con i tempi esso è oggi "in rete" grazie alla realizzazione e al costante aggiornamento del sito web *oulipo.net* poiché «il suo vero segreto risiede da sempre nella sua assoluta trasparenza».[72]

*Metodo di lavoro collettivo*. Per entrambi i gruppi una prassi di lavoro comune è un punto saldo. Nel caso dell'Oulipo emblematico è il seguente scambio di battute:

> «QUENEAU: Non dobbiamo mai dimenticare che i lavori dell'OuLiPo sono comuni. Il loro risultato è il bene di tutti. Nessuno può disporne senza il consenso di tutti.
> LE LIONNAIS: Lo stesso vale per qualsiasi opinione. Perché tutto riguarda l'Oulipo, compresa la conferenza di Evian e la fabbricazione della *gibelotte*.
> LATIS: Benissimo».[73]

Tuttavia, parallelamente, ciascun *oulipiens* è libero di dedicarsi alla propria carriera così come di promuovere attività individuali e pubblicare opere composte secondo i principi compositivi della *LiPo* o meno. In questo modo il collettivismo non esclude l'individualismo. Testimonianza di ciò sono ad esempio le *Bibliographies des membres de l'Oulipo* presenti nell'*Atlas* [27, p. 409-426]. In esse, per ciascun membro, sono elencate tanto le opere oulipiane (o relative all'Oulipo) quanto una selezione di opere non oulipiane o solo parzialmente oulipiane. Infine, nonostante tra i membri dell'Oulipo ci sia un certo *senso di filiazione* tant'è che, per i membri delle generazioni successive a quella dei fondatori, viene spesso indicato da chi sono stati introdotti, vige tuttavia l'autonomia intellettuale e operativa ben esplicitata da Roubaud quando afferma che «Raymond Queneau è il mio maestro, ma sono io che decido e so in cosa, come e in che misura»[74].

Diversa è invece l'accezione di lavoro collettivo professata da Bourbaki. Il processo di revisione tra pari a cui è sottoposto, per esempio, ciascuno dei volumi del *Traité* è estremante rigido e il testo può essere pubblicato solo dopo essere stato rimaneggiato praticamente da ciascun bourbakista e approvato all'unanimità. Questo finisce per rendere indistinguibile il contributo dell'uno o dell'altro: il singolo si diluisce nella collettività. Queneau descrive questo processo creativo col suo abituale e lucido sarcasmo: «Come lavorano i bourbakisti? Si dice che [...] uno di loro sia incaricato della stesura di un fascicolo. Questa prima versione viene inviata ai vari membri di Bourbaki che poi si riuniscono al Congresso; esaminano questa prima bozza, la criticano: in generale non ne rimane nulla. E così via, fino a quando questi nulla integrati [tra loro] danno un insieme non vuoto di teoremi, proposizioni, annotazioni e, possibilmente, assiomi tutti nulli eccetto un numero finito».[75]

---

[69] *Notice sur la vie et l'oeuvre de Nicolas Bourbaki* http://archives-bourbaki.ahp-numerique.fr/files/original/e460d9db53e0e52cecf98c9049b1bc69.pdf (Accessed 22 March 2021).
[70] Italo Calvino, *Perec, gnomo e cabalista*, La Repubblica, 6 marzo 1982, p.18.
[71] https://oulipo.net/fr/jeudis (Accessed 22 March 2021).
[72] «L'Ouvroir de Littérature Potentielle a longtemps passé pour une organization secrète. Mais son véritable secret a toujours résidé dans son absolue transparence» [27, p. 407].
[73] «QUENEAU: Nous ne devons jamais oublier que les travaux de l'OuLiPo sont *communs*. Leur résultat est le bien de tous. Nul ne peut en disposer sans l'accord de tous. Aggiungere chiarimento
LE LIONNAIS: Il en est de même de toute opinion. Car tout concerne l'Oulipo, y compris la conférence d'Evian et la fabrication de la gibelotte.
LATIS: Ben alors» [9, p. 71-72].
[74] «Raymond Queneau est mon maître, mais c'est moi qui décide et sais en quoi, comment et jusqu'où» [*La boucle*. Seuil, Paris (1993), p.270].

Certamente, è comune ai due gruppi la volontà di coltivare la collegialità delle attività e rafforzare il senso di appartenenza al gruppo. Questo, nel caso di Bourbaki, serve «molto concretamente, [anche a] eliminare eventuali diritti d'autore in un fondo comune [mentre] la preoccupazione dell'Oulipo è [soprattutto] quella di preservare la dimensione innovativa dell'impresa collettiva, senza gravarla fin dall'inizio del peso della notorietà del tale e del tal altro» [25].

*Organi di stampa*. Un ulteriore parallelismo tra i due gruppi può essere fatto sul piano dei mezzi di diffusione delle opere.

All'interno di Bourbaki le comunicazioni sono state veicolate dal 1940 al 1953 tramite i fascicoli del *Journal de Bourbaki* detto *La Tribu*. Dapprincipio i fascicoli erano ovviamente a esclusivo uso interno ma oggi l'archivio è stato desecretato e alcuni volumi sono disponibili addirittura sul web. La struttura e il contenuto di ciascun volumetto è in buona sostanza quella di un verbale di un incontro: con fare oulipiano potremmo dire che si tratti di un ulteriore caso di plagio per anticipazione dei *Comptes rendus* delle riunioni degli *oulipiens* (in [9] sono raccolti quelli dei primi 3 anni di vita dell'*Ouvroir*).

La *Bibliothèque oulipienne* è invece una collana di fascicoli che raccoglie i lavori degli *oulipiens* sia individuali che collettivi. Ad oggi conta più di 200 fascicoli di ciascuno dei quali vengono stampati solo 150 esemplari numerati (più 50 riservati ai membri) ed è già previsto un fascicolo n. 666 a firma collettiva intitolato *Diable!* naturalmente. I fascicoli sono stati periodicamente riuniti in volumi e pubblicati da Seghers fino al 1990, poi da Castor Astral. Memorabili il n. 1 *Ulcérations* di Perec e il n. 2 *La princesse Hoppy ou Le conte du Labrador* di Roubaud entrambi del 1974, così come il n. 6 *Piccolo Sillabario Illustrato* di Calvino del 1978, giusto per citarne alcuni[76].

L'opera con cui, a partire dal 1939 e fino al 1998, Bourbaki si svela al pubblico sono gli 11 libri degli *Éléments de mathématique* sui quali molto è stato scritto e pertanto si rimanda a bibliografia specifica (per esempio [6] e [7]).

Sul fronte dell'Oulipo le due opere collettive di maggiore respiro e ampia diffusione sono sicuramente *La Littérature potentielle: créations, re-créations, récréations* del 1973 [28] e, sette anni dopo, l'*Atlas de littérature potentielle* [27]. Del 2009 è invece l'*Anthologie de l'OuLiPo* [26] che ripropone una selezione dei lavori oulipiani e qualche inedito.

A un lettore attento delle opere dei due gruppi non sfuggiranno alcuni rimandi: un esempio illustre è dovuto a Perec che col suo *La Vie mode d'emploi* sembra citare l'avvertenza che apre ogni fascicolo del *Traité* intitolata appunto *Mode d'emploi de ce traité*.

*Universalità*. «Bourbaki e Oulipo hanno in comune la rivendicazione di una potenziale universalità».[77] Infatti, così come la matematica di Bourbaki chiaramente non si limita a un solo luogo o a una sola lingua, allo stesso modo la pratica della scrittura "vincolata" è concepibile in tutte le lingue anche se alcune *contrainte* possono non essere applicabili a tutti gli idiomi. È possibile farsi un'idea di questo aspetto leggendo, ad esempio, le riflessioni di Umberto Eco sulle sfide affrontate e sulle decisioni operate nel tradurre gli *Esercizi di stile* di Queneau.[78]

---

[75] «Comment travaillent les bourbakistes? On dit que [...] l'un d'entre eux soit chargé de la rédaction d'un fascicule. Cette première mouture est envoyée aux différentes Bourbaki qui se réunissent ensuite en Congrès; on examine cette première rédaction, on la critique: en général il n'en reste rien. Et ainsi de suite, jusqu'à ce que ces riens intégrés donnent un ensemble non vide de théorèmes, propositions, scholies et, éventuellement, axiomes tous nuls à l'exception d'un nombre fini» [30, p.7].
[76] Per un elenco completo si può fare riferimento al sito dell'Oulipo, in particolare alla pagina: https://www.oulipo.net/fr/bo (Accessed 27 April 2021)
[77] «The mathematics of Bourbaki, even if some of their aspects were attributed (whether to praise or belittle them) to a certain "French genius", are clearly not limited to a single land or tongue. The practice of writing by constraints is conceivable - even if certain constraints cannot be generalized everywhere - in all languages. Bourbaki and Oulipo have in common the claim to a potential universality» [34, p. 128].

In questo contesto vale la pena osservare inoltre come l'Oulipo abbia trasceso i confini nazionali francesi e, per esempio, in Italia abbia ispirato nel 1990 la nascita dell'*Opificio di Letteratura Potenziale* (OpLePo) i cui fondatori furono Raffaele Aragona, Ruggero Campagnoli e Domenico D'Oria mentre, tra i suoi membri, si annoverano tra gli altri Edoardo Sanguineti, Paolo Albani, Piergiorgio Odifreddi. L'Oulipo ha anche saputo travalicare i confini propriamente detti della Letteratura. Sul modello dell'Oulipo infatti Le Lionnais immagina altri *Laboratori* [21]*, Ou-*x*-po*[79], per la musica (*x=mu*), la pittura (*x=pein*), il cinema (*x=ciné*), il teatro tragicomico (*x=tra*), la cucina (*x=cui*), i fumetti (*x=ba*) e, il 23 agosto del 1973, fonda in prima persona l'*Ouvroir de Littérature Policière Potentielle*, noto come Oulipopo [42, p. 398].

*Reclutamento.* In entrambi i gruppi l'ingresso dei nuovi membri avviene all'unanimità e per cooptazione. Quest'ultima in Bourbaki è impegnativa ed esigente, una sorta di rito iniziatico. I giovani matematici ritenuti più promettenti vengono invitati ai Congressi di Bourbaki (tre ogni anno) e ricevono il titolo temporaneo di "cavie" (fr. *cobaye*). Partecipano alle sessioni dove intervengono, fanno e ricevono domande, ma solo coloro i quali soddisfano gli elevati standard bourbakisti sono infine introdotti nel gruppo. La cooptazione nell'*Ouvroir* è di contro più *soft*: non si può chiedere di essere ammessi ma si viene individuati, proposti da un membro, "osservati" e infine cooptati se il giudizio è unanime. Ad oggi si contano 38 affiliati inclusi i membri morti «excusés pour cause de décès» secondo una espressione coniata da Le Lionnais stesso. L'ultimo *oulipien* reclutato nel 2014 è lo scrittore argentino Eduardo Berti; prima di lui, nel 2012, il fumettista Étienne Lécroart membro dell'*Oubapo*.

L'attitudine alla classificazione e la passione per la combinatoria proprie del gruppo hanno portato Roubaud a individuare quattro differenti tipi di membri che, facendo uso degli operatori propri della teoria degli insiemi (con un'intersezione non simmetrica), ci piace denotare come:

1. **L\M:** letterati (scrittori di prosa, poeti, critici e giornalisti) che *non sono* matematici;
2. **M\L:** matematici che *non sono* letterati;
3. **L∩M**: letterati *e* matematici (come Roubaud stesso e Queneau e, tra i membri dell'ultima generazione Salon);
4. **M∩L**: matematici *e* letterati.

Egli specifica che, nei casi 3 e 4, la *e* non va intesa come il connettore logico *and* ma piuttosto come *and then*. Ciò fa sì che «l'Oulipo non abbia un ruolo esclusivamente letterario»[80] tanto più che i matematici svolgono un ruolo tutt'altro che marginale nella vita del gruppo. Tra questi Berge, co-fondatore ed esperto di teoria dei grafi, autore del romanzo poliziesco *Qui a tué le Duc de Densmore?* [11] il cui intrigo può essere risolto solo da un teorico dei grafi grazie all'applicazione di un teorema del matematico ungherese Elias M. Hagos. Ne *La Princesse Aztèque* [10][81] espone invece le regole per la composizione di un *sonetto di lunghezza variabile* proponendo le trasformazioni geometriche utili a riscrivere un sonetto di 14 versi e 12 sillabe in uno di 15 versi e 12 sillabe. Ciò è realizzabile solo agendo sulla già "vincolata" struttura del sonetto con un numero non indifferente di *contrainte* aggiuntive quali restrizioni prosodiche, sintattiche, semantiche, restrizioni sulla disposizione delle righe e altre ancora.

Tra gli esponenti dell'insieme 4 ricordiamo anche Paul Braffort logico ed esperto d'informatica che vanta l'onore di essere stato il primo cooptato[82], Pierre Rosenstiehl teorico dei labirinti. Scorrendo i nomi degli invitati alle riunioni dell'*Ouvroir* si incontrano, per esempio, il bourbakista René Thom noto per la sua Teoria delle catastrofi e Jean Ferry ex-patafisico ed ex-cavia di Bourbaki, e tutto ciò manifesta chiaramente l'interesse del gruppo ad annettere quanti più matematici possibile.

---

[78] «Si trattava, in conclusione, di decidere cosa significasse, per un libro del genere, essere fedeli. Ciò che era chiaro è che non voleva dire essere letterali. Diciamo che Queneau ha inventato un gioco e ne ha esplicitato le regole nel corso di una partita, splendidamente giocata nel 1947. Fedeltà significava capire le regole del gioco, rispettarle, e poi giocare una nuova partita con lo stesso numero di mosse» [31, *Introduzione* p. XIX].

[79] Per un elenco delle attività degli Ou-*x*-Po si veda il sito: http://www.fatrazie.com/jeux-de-mots/ouxpo (Accessed 22 April 2021).

[80] «This means that the Oulipo *does not have* an exclusively literary role» [34, p. 126].

[81] Anche in [26, p. 763].

[82] [9, p. 47].

Oltre alla capacità di rinnovarsi e reiventarsi – secondo Montémont: «l'Oulipo, dopo un'eclissi negli anni '90, deve la sua sopravvivenza in gran parte all'energia che Marcel Bénabou ha profuso per migliorarne la visibilità pubblica ed editoriale»[83] – la longevità dell'Oulipo è anche da attribuire - seppur in minima parte - all'impossibilità di lasciare il gruppo. Per volere di Queneau, segnato negativamente dalla sua esperienza all'interno del gruppo surrealista[84], gli *oulipiens* si sono dotati di una specifica procedura secondo cui nessun membro può cessare di appartenere all'Oulipo neanche da morto, se non suicidandosi davanti a un ufficiale giudiziario e dichiarando che il gesto ha proprio «lo scopo di liberarlo dall'Oulipo e ripristinare la sua libertà di manovra per il resto dell'eternità».[85]

Di contro, ad ogni bourbakista è imposto di ritirarsi dal gruppo una volta compiuti i cinquant'anni a causa della (falsa) credenza secondo cui un matematico è veramente creativo e produttivo solo da giovane.

**Conclusione**

Oggi l'Oulipo è più vivo e prolifico che mai mentre Bourbaki ha cessato di svolgere un ruolo da protagonista nel panorama della matematica mondiale, pur lasciando una preziosissima eredità culturale. Citando Queneau infatti potremmo dire che «la matematica procede a grandi passi e Bourbaki rappresenta gli stivali delle sette leghe che chiunque deve calzare se vuole raggiungerlo (il treno)».[86]

Montémont riferisce l'opinione di Roubaud secondo cui «l'Oulipo non è solo un'emanazione derivata da Bourbaki: ha saputo premunirsi, molto intelligentemente, contro l'isolamento intellettuale, il settarismo o il rifiuto di qualsiasi applicazione concreta, altrettanti errori che hanno portato il matematico policefalo al suo declino».[87]

Del resto, la natura del compito che Bourbaki si era dato, ovvero ricostruire e unificare l'intera Matematica sulla sola base del metodo assiomatico e di poche "strutture madri" fondamentali e delle loro combinazioni, era intrinsecamente tale da condurre alla dissoluzione del gruppo una volta portato a termine il lavoro. Ci piace però pensare, sulla scorta di Queneau, che Bourbaki sia "morto" senza mai invecchiare: «è per forza invecchiato il tuo matematico immaginario, deve essere rimasto indietro. Beh! no, Bourbaki non è invecchiato perché non può invecchiare».[88] È dunque tramontato solo il mito di Bourbaki [2] mentre è destinato all'eternità ciò egli ha rappresentato sia all'interno della comunità matematica che al di fuori di essa, dal momento che «come simbolo […] è stato, per più di 30 anni, abbastanza potente da servire a molti scopi diversi in tutte le discipline. Osservando i vari ruoli che ha svolto in diversi tipi di trattazioni, il suo crescente impatto su una porzione di matematici, strutturalisti e scrittori allo stesso modo, e quindi la sua autorità in declino, possiamo studiare il modo in cui diversi flussi culturali si mescolavano in un nodo chiamato Bourbaki».[89] David Aubin considera Bourbaki come un *connettore culturale* ovvero come un riferimento autorevole, più o meno esplicito, attinto da altre discipline e usato quando si vuole argomentare su un tema o conferire maggiore legittimità ai propri metodi e alle proprie idee.[90] A tal proposito, ci sembra appropriato concludere con le parole di Le Lionnais: «il movimento bourbakista è un fatto; […] gli si deve essere riconoscenti. Il suo tentativo non costituisce che uno dei modi di continuare la matematica dei nostri giorni».[91] In questo senso, per un autore invocare Bourbaki equivale a importare nel proprio ambito culturale interi insiemi di significati e pratiche. Ciò rende evidenti le intersezioni culturali tra matematica bourbakista, Strutturalismo e letteratura potenziale.

---

[83] «L'Oulipo, lui, après une éclipse dans les années 90, a dû sa survie, en grande partie, à l'énergie qu'a mise Marcel Bénabou pour en améliorer la visibilité publique et éditoriale» [25, p. 10].

[84] Queneau acconsentì alla fondazione a condizione di «scartare radicalmente ogni attività del gruppo che possa generare fulminazioni, scomuniche e ogni forma di terrore» («écarter de manière radicale toute activité de groupe pouvant engendrer fulminations, excommunications et toute forme de terreur») [23, p. 77].

[85] «his suicide was intended to release him from the Oulipo and restore his freedom of manoeuvre for the rest of eternity» [34, p. 126].

[86] «la mathématique continue à advancer à grand train, et Bourbaki est les bottes de sept lieues que quiconque doit chausser s'il veut le rattraper (le train)» [30, p. 9].

[87] «l'Oulipo n'est pas seulement une émanation dérivée de Bourbaki: il a su se prémunir, fort intelligemment, contre l'isolement intellectuel, le sectarisme ou le refus de toute application concrète, autant de travers qui ont conduit le mathématicien polycéphale à son déclin» [24]. L'appellativo di "matematico policefalo" è stato attribuito a Bourbaki da Carl Boyer in *Storia della matematica*, ISEDI, Torino (1976), p. 718.

[88] «Il a nécessairement vieilli votre fictif mathématicien, il doit avoir pris du retard. Eh bien! non, Bourbaki n'a pas vieilli parce qu'il ne peut pas vieillir» [30, p. 7].

[89] «As a symbol […] he was, for more than 30 years, powerful enough to serve many different purposes across disciplines. By looking at the various roles he played in several types of discourse, his rising impact among a portion of mathematicians, structuralists, and writers alike, and then his declining authority, we can study the way in which different cultural streams mingled at a node called Bourbaki. We will be able to put some flesh on the cliche that sometimes, somehow, ideas are "in the air"» [2, p. 297].

---

[90] «I defined cultural connectors […], as more or less explicit references from other disciplines, used by actors when they attempt to argue for a point or when they want to increase the legitimacy of their methods and ideas» [1].
[91] «Le mouvement bourbakiste est un fait […] on doit lui être reconnaissant. Sa tentative ne constitue que l'un des moyens de continuer les mathématiques de nos jours» [18, p.10].